\documentclass[11pt]{amsart}
\usepackage{amssymb,amsthm,latexsym,amsfonts}
\usepackage[latin1]{inputenc}
\usepackage[T1]{fontenc}

\def\n{{\nabla_{\kern-4pt\scriptscriptstyle H}}}

\def\TN{\tilde{N}}

\defø{{\scriptscriptstyle\otimes}}

\def\F{{\mathfrak F}}

\def\P{{\mathbf P}}
\def\N{{\mathbb N}}
\def\I{{\mathbf 1}}

\def\R{{\mathbb R}}

\def\integ{\int_0^T\kern-7pt\int_E}

\def\({\Bigl(}
\def\){\Bigr)}
\def\n{\ \nu(ds,dz)}
\defø{\ \omega(ds,dz)}
\defØ{\Omega}

\makeatletter
\newenvironment{example*}[1][\examplename]{\par
  \normalfont
  \topsep6\p@\@plus6\p@ \trivlist
  \item[\hskip\labelsep\itshape
    Example #1\@addpunct{.}]\ignorespaces
}{%
  $\sqsupset$\endtrivlist
}

\newcommand{\examplename}{Example}
\makeatother

\newtheorem{theo}{Theorem}

\newtheorem{remar}{Remark}

\newtheorem{prop}{Proposition}

\newtheorem{propo}{Proposition}
\newtheorem{lemma}{Lemma}
\newtheorem{hyp}{Hypothesis}

\begin{document}

\title{Girsanov Theorem for Filtered Poisson Processes}
\author{L. Decreusefond}
\address{
Département Informatique et Réseaux\\
 Ecole Nationale Supérieure des Télécommunications\\
46, rue Barrault, 75634 Paris Cedex 13, France 
}

\author{N. Savy}
\address{
Institut de Recherche Mathématiques de Rennes\\
Université  de Rennes 1, 35042 Rennes Cedex, France
}

\begin{abstract}

Shot-noise and fractional Poisson processes are instances of filtered
Poisson processes. We here prove  Girsanov theorem for this kind of
processes and give an application to an estimate problem. 
\end{abstract}

\maketitle
\section{Introduction}
At the beginning of the theories of mathematical finance and
performance analysis of telecommunications networks, asset prices or
input traffics were modeled by Markovian diffusion processes. It is
by now well known that these processes don't capture some of the
essential features observed in the real datas. Namely, diffusion
processes do no reflect long-range dependance exhibited in many traces
\cite{kuhn02,leland94}. Alternative models are currently being
considered such as L\'evy processes, fractional Brownian motion, shot
noise processes and different combinations of them. For instance, in
\cite{kuhn02}, it is suggested that one could represent the evolution
of an asset price by a process $S_t=\exp(B_t+N_t),$ where $B$ is a
standard Brownian motion and $N$ is a shot noise process (for
discussions about the validity of this model we refer to
\cite{kuhn02}). Motivated by this work, we here investigate some of
the basic properties of processes which encompass shot-noise models.
 Given a
marked Poisson process $N$ and a deterministic kernel $K$, we define 
filtered Poisson process $N^K,$ by
\begin{align*}
N^K_t  &=\int_0^t K(t,s) dN_s \\
        &=\int_0^t \int_E z K(t,s) \mu(ds,dz). 
\end{align*} 
If $K$ is a convolution kernel, i.e., $K(t,s)=k(t-s)$ then $N^K$ is
usually called a shot-noise process.
Shot-noise processes are used in numerous fields of applications, e.g. electronics, hydrology, climatology,
telecommunications (see \cite{parzen}, \cite{snyder}, \cite{yh01} and 
references therein), insurance (see \cite{dj} and references therein)
and finance (see \cite{kuhn02} and  \cite{s96} for a review on this topic).
It follows that filtered Poisson processes can, in the same fields,
model more general phenomena so that we believe, it is an interesting
class of processes to consider. This paper is organized as follows. In the next section, we give some
preliminaries on point processes. In section 3, we prove the Girsanov
theorem for filtered Poisson processes and in section 4, we apply this
result to an estimate problem.

\section{Preliminaries}
Let $\R^d$ equipped with its borelian $\sigma$-algebra and
$(Ø,\F,(\F_t)_{t \in \R^+},\P)$ a filtered probability space. 
We consider that we are given an $\F$-adapted, marked Point process
$\mu$ of intensity measure $\nu$ of the form (for details, see \cite{jacod})
\begin{displaymath}
\nu(\omega,ds,dz)=\lambda(s) ds \,  \eta(dz).
\end{displaymath}
This means that there exist  two sequences of random variables $(T_n)_{n \in \N^*}$
and $(Z_n)_{n \in \N^*}$ such that : 
$(T_n)_{n \in \N^*}$ is  a strictly increasing sequence of $\R^+$
(jump times) verifying $\lim_{n \to \infty}T_n=+ \infty$ and 
$(Z_n)_{n \in \N^*}$ is a sequence of $\R^d$ (the marks associated
with the jump times), such that 
\begin{displaymath}
\mu_{t,z}=\sum_{n \in \N^*} \delta_{(T_n,Z_n)}(t,z).
\end{displaymath}
We associate with $\mu$ the Marked Poisson process $N$ defined as:
\begin{displaymath}
N_{t}=\sum_{n \in \N^*} Z_n\I_{[T_n \leq t]}.
\end{displaymath}

\begin{hyp} \label{H-HK1}
We assume now that our deterministic kernel satisfies a few regularity
assumptions :
\begin{itemize}
\item $K$ is triangular, i.e., $K(t,s)=0$ for $s>t.$ 
\item $K$ is sufficiently integrable in the sense that $(s,z) \to zK(t,s)$ is in ${\mathcal{L}}^{2}(\nu)$ for any
$t> 0$.
\item $K$ is continuously differentiable with respect to $t$ and $s$
  in $\{(t,s):\, 0<s<t\}.$
\end{itemize}
\end{hyp}
We now define a filtered Poisson process by 
\begin{equation}
  \label{eq:1}
  N^K_t=\int_0^t \int_E z\, K(t,s) \mu(ds,dz).
\end{equation}
 The first hypothesis thus appears as necessary if we want to have a
 non- anticipative process $N^K.$ The other two are only technical
 assumptions. 
The next proposition is immediate.
 \begin{prop}
   The sample-paths of $N^K$ are c\`adl\`ag iff $K(t,t)$ is
finite for any $t$ and that they are continuous iff $K(t,t)=0$
for any $t.$ 
 \end{prop}

\section{Girsanov Theorem}
We hereafter assume that $K$ is degenerate on the diagonal, that is to
say $K(t,t)=0$ for any $t.$ It has thus a sense to look at the
perturbations of the sample-paths which induce an absolutely
continuous change of probability measure.
\begin{theo} \label{T-GIR}
Let $\{\tilde{N}_t:t \geq 0\}$
the compensated Marked Poisson Process associated with $N$ defined by
\begin{displaymath}
  \tilde{N}^K_t=\int_0^t \int_E z\, K(t,s) (\mu(ds,dz)-\nu(ds,dz)).
\end{displaymath}

Consider a function $h\in \mathcal{L}^1(\nu)$ and introduce the
probability measure $\P_h$ which is absolutely continuous with respect
to $\P$ of Radon-Nikodym  density given by: 
\begin{align}
Z_t   
  &=\left.{\frac{d\P_h}{d\P}}\right|_{{\mathfrak F}_t} \notag\\
  &= {\mathcal{ED}}\(\int_0^.\!\!\int_E h(s,z) (\mu-\nu)(ds,\ 
dz)\)_t \label{E-DOLEANS}
\end{align}
where $\mathcal{ED}(X)_t$ indicates the Doléans-Dade exponential of the
process $X$ at time $t$.

Let $N^{h,K}$ be defined by
\begin{displaymath}
  N^{h,K}_t=\int_0^t \int_E z\, K(t,s)
  (\mu(ds,dz)-\nu(ds,dz))-\int_0^t \int_E z\, K(t,s)h(s,z)\, \nu(ds,dz).
\end{displaymath}
Then the following equality in law holds:
\begin{displaymath}\mathcal{L}(\tilde{N}^K,\P_h)=\mathcal{L}(\tilde{N}^{h,K},\P)\end{displaymath}
\end{theo}
\begin{proof}
Let $n \in \N^*$ and $(t_1, \dots, t_n) \in [0,T]^n$ fixed. When these
instants are fixed, the processes
\begin{displaymath}
M^i_r:r \to \int_0^r \int_E z K(t_i,s) \mu(ds,dz)
\end{displaymath}
are marked Point processes whose compensators are:
\begin{displaymath}
\Phi^i_r:r \to \int_0^r  \int_E z K(t_i,s) \nu(ds,dz)
\end{displaymath}
Hence, $\left(M^i-\Phi^i\right)_{i=1}^n$ is a $\P$-martingale and the
Girsanov theorem (see \cite{jacod}) makes certain that:
\begin{displaymath}
\left( M^{i,h}:r \to M^i_r-\Phi^i_r - \int_0^r \frac{1}{Z_{s^-}}
d \left<M^i-\Phi^i,Z\right>_s \right)_{i=1}^n
\end{displaymath}
is a  $\P_h$-martingale.
$Z$ is defined by \eqref{E-DOLEANS} so $Z$ is solution of the equation:
\begin{displaymath}
R_t=1+\int_0^t \int_E R_{s^-}[h(s,z)-1](\mu - \nu)(ds,dz)
\end{displaymath}
Hence, we can write:
\begin{align}
(A)    &=\left<M^i-\Phi^i,Z\right>_s \notag\\
        &=\left<\int_0^.\int_E z K(t_i,s) (\mu - \nu)(ds,dz),1+\int_0^. \int_E 
Z_{s^-}[h(s,z)-1] (\mu - \nu)(ds,dz) \right>_s \notag\\
        &=\left<\int_0^. \int_E z K(t_i,s) (\mu - \nu)(ds,dz),\int_0^. \int_E  
Z_{s^-}[h(s,z)-1] (\mu - \nu)(ds,dz) \right>_s \notag\\
        &=\int_0^s z  K(t_i,s) Z_{s^-}[h(s,z)-1] 
d\left< \int_0^. \int_E (\mu - \nu)(ds,dz) ; 
\int_0^. \int_E  (\mu - \nu)(ds,dz) \right>_s \notag\\
        &=\int_0^s \int_E z K(t_i,s) Z_{s^-}[h(s,z)-1] \nu(ds,dz)  \notag
\end{align}         
We can conclude that:
\begin{align}
M^{i,h}_r
        &=M^i_r - \Phi^i_r -  \int_0^r \frac{1}{Z_{s^-}}d
\left<M^i - \Phi^i,Z\right>_s \notag\\
        &=M^i_r - \Phi^i_r - \int_0^r \frac{1}{Z_{s^-}}\int_E z K(t_i,s)
Z_{s^-}[h(s,z)-1] \nu(ds,dz)
dz \notag\\
        &=M^i_r - \Phi^i_r - \int_0^r \int_E z K(t_i,s) [h(s,z)-1] \nu(ds,dz)  \notag\\
        &=\int_0^r \int_E z K(t_i,s) \mu(ds,dz)- \Phi^i_r - \int_0^r
\int_E z K(t_i,s) h(s,z)
\nu(ds,dz) + \Phi^i_r  \notag\\
&=\int_0^r \int_E z K(t_i,s) \mu(ds,dz) -  \int_0^r \int_E z K(t_i,s) h(s,z)
\nu(ds,dz)      \label{E-Julie}
\end{align}
thus,
\begin{displaymath}
r \to \int_0^r \int_E z K(t_i,s) \mu(ds,dz) - \int_0^r \int_E
z  K(t_i,s)
h(s,z) \nu(ds,dz)
\end{displaymath}
is a $\P_h$-martingale and 
\begin{equation} \label{E-Louis}
r \to \int_0^r \int_E z K(t_i,s) h(s,z) \nu(ds,dz)
\end{equation}
is the $\P_h$-compensator of $M^i$.\\
>From \eqref{E-Julie}, $M^{i,h}$ is written:
\begin{displaymath}
M^{i,h}:r \to \int_0^r \int_E z K(t_i,s) (\mu(ds,dz) - h(s,z)
\nu(ds,dz))
\end{displaymath}
So its compensator under $\P$ is \eqref{E-Louis}.\\
It follows (see \cite{jacod}) that: 
\begin{displaymath}
\mathcal{L} \left( \left[r \to M^i_r \right]_{i=1}^n ;
\P_h \right) =\mathcal{L} \left( \left[r \to M^{i,h}_r \right]_{i=1}^n ;
\P \right)
\end{displaymath}
Now, taking $r = \underset{1 \leq i \leq n}{Sup} \, t_i$, the kernel being triangular, the result holds.
\end{proof}

\section{An application}
Suppose now that we are given to observe the paths of a perturbed
filtered Poisson process:
\begin{displaymath}
  X^{\theta}_t=\tilde{N}^{K}_t-\theta t \text{ with } \theta >0.
\end{displaymath}
We now investigate the problem of the estimate of $\theta.$ Such a
problem has been thoroughly investigated when a standard Brownian motion
is put in place of $N$ (i.e., a fractional Brownian motion is
substituted to  $\tilde{N}^K$) \cite{du,MR1897916,MR2001k:93131}. We
here  need to
assume that $\nu$ and $K$ are such that there exists one and only one
function $\phi$ such that:
\begin{equation*}
  \int_0^t \int_E z K(t,s) \phi(s) \nu(ds,dz)= t.
\end{equation*} 
We then  introduce the likelihood function:
\begin{equation*}
Z^{\theta}_t= \left.{\frac{d\P_\theta}{d\P}}\right|_{{\mathfrak
    F}_t}  =\mathcal{ED} \left(\int_0^. \int_E \theta \phi(s)
(\mu-\nu)(ds,dz) \right)_t.
\end{equation*}
According to the previous theorem, the processes $X^{\theta}$ under $\P$ and $\TN^K$ under
$\P_{\theta}$ have the same law thus an estimate of  $\theta$ is  given by:
\begin{equation*}
\hat{\theta}_t=\underset{\theta \in [0,1]}{Argmax} \, Z^{\theta}_t.
\end{equation*}
Unfortunately, the exact expression of the Dol\'eans-Dade exponential
is so intricate in case of jump processes that we have to find another
expression of $Z^\theta$ more suitable for computations $\hat{\theta}_t.$
\begin{lemma}
If $\phi \in \mathcal{L}^1(\nu)$ and  $\ln(1+\theta \phi) \in
\mathcal{L}^1(\nu)$ then:
\begin{equation*}
\left\{
\mathcal{ED} \left(\int_0^. \int_E \theta \phi(s)
(\mu-\nu)(ds,dz) \right)_t : 0 \leq t \leq T \right\}
=
\left\{
\exp(Y^{\theta}_t)  : 0 \leq t \leq T \right\}
\end{equation*}
with for any  $0 \leq t \leq T$:
\begin{equation*}
Y^{\theta}_t=\int_0^t \int_E \ln(1+\theta \phi(s)) \mu(ds,dz)-\int_0^t \int_E
\theta \phi(s) \nu(ds,dz)
\end{equation*}
\end{lemma}

\begin{proof}
Apply Itô formula to the function $x
\to \exp(x)$ with respect to the jump process $Y^{\theta}$. One gets:
\begin{multline*}
\exp(Y^{\theta}_t) - 1\\
        \begin{aligned}
        = &\int_0^t \exp(Y^{\theta}_{s^-}) dY^{\theta}_s
        +\sum_{s \leq t} [\exp(Y^{\theta}_s) -
        \exp(Y^{\theta}_{s^-}) - \exp(Y^{\theta}_{s^-})(Y^{\theta}_s -
        Y^{\theta}_{s^-})] \\
        = &\int_0^t \exp(Y^{\theta}_{s^-}) dY^{\theta}_s
        +\int_0^t \int_E [\exp(Y^{\theta}_{s^-} + \ln(1+\theta \phi(s)))-\exp(Y^{\theta}_{s^-}) \\
        & \qquad - \exp(Y^{\theta}_{s^-})(\ln(1+\theta \phi(s))) \mu(ds,dz)] \\
        = &\int_0^t \exp(Y^{\theta}_{s^-}) dY^{\theta}_s
        +\int_0^t \int_E [\exp(Y^{\theta}_{s^-} + \ln(1+\theta \phi(s))) \\
        & \qquad - \exp(Y^{\theta}_{s^-})(1 + \ln(1+\theta \phi(s))) \mu(ds,dz)] \\
        \end{aligned}
\end{multline*}
But we have:
\begin{equation*}
dY^{\theta}_s= \int_E \ln(1+\theta \phi(s)) \mu(ds,dz)- \int_E
\theta \phi(s) \nu(ds,dz)
\end{equation*}
hence:
\begin{multline*}
\exp(Y^{\theta}_t) - 1\\
        \begin{aligned}
        = &\int_0^t \int_E \exp(Y^{\theta}_{s^-}) [ \ln(1+\theta \phi(s)) \mu(ds,dz)- 
        \theta \phi(s) \nu(ds,dz)] \\
        & \qquad +\int_0^t \int_E [\exp(Y^{\theta}_{s^-} + \ln(1+\theta \phi(s))) \\
        & \qquad - \exp(Y^{\theta}_{s^-})(1 + \ln(1+\theta \phi(s)))]
        \mu(ds,dz) \\
        = &\int_0^t \int_E \exp(Y^{\theta}_{s^-}) [ \ln(1+\theta \phi(s))-(1 +
\ln(1+\theta \phi(s)))\\
        & \qquad +\exp( \ln(1+\theta \phi(s)))] \mu(ds,dz) \\
        & \qquad -\int_0^t \int_E \exp(Y^{\theta}_{s^-}) \theta \phi(s)
        \nu(ds,dz) \\
        = &\int_0^t \int_E \exp(Y^{\theta}_{s^-}) \theta \phi(s) \mu(ds,dz)
        -\int_0^t \int_E \exp(Y^{\theta}_{s^-}) \theta \phi(s)
        \nu(ds,dz) \\
        = &\int_0^t \int_E \exp(Y^{\theta}_{s^-}) \theta \phi(s)
        (\mu-\nu)(ds,dz)\\
                \end{aligned}
\end{multline*}
It follows that the process  $\exp(Y^{\theta})$ is a solution to the Stochastic
Differential Equation which defines the Dol\'eans-Dade exponential. By
uniqueness of the solution, the result holds.
\end{proof}
\begin{propo}
Suppose
$\phi \in \mathcal{L}^1(\nu)$ and $\ln(1+\theta\phi) \in
\mathcal{L}^1(\nu)$. Then,  for any $t \in [0,T],$ $\hat{\theta}_t$
exists and 
is the unique positive solution of 
\begin{equation}\label{eq:5}
  \sum_{T_j\le t} \frac{\phi(T_j)}{1+\hat{\theta}_t\phi(T_j)}=\int_0^t \phi(s)\lambda(s) ds.
\end{equation}
\end{propo}
\begin{proof}
Fix $t \in [0,T]$.
The likelihood function reaches its maximum when the function:
\begin{equation*}
f:\theta \to\int_0^t \int_E \ln(1+\theta \phi(s)) \mu(ds,dz)-\int_0^t
\int_E
\theta \phi(s) \nu(ds,dz)
\end{equation*}
reaches its maximum. Noticing that:
\begin{align*}
f'(\theta)     &=\int_0^t \int_E \frac{\phi(s)}{1+\theta
                \phi(s)} \mu(ds,dz)-
                \int_0^t \int_E  \phi(s) \nu(ds,dz) \\
f''(\theta)    &=\int_0^t \int_E \frac{-\phi(s)^2}{[1+\theta
                \phi(s)]^2} \mu(ds,dz),
\end{align*}
wee see that the function $f$ is concave  hence  that it admits a
unique maximum and that this maximum $\hat{\theta}_t$ satisfies \eqref{eq:5}.
\end{proof}
\begin{remar}
If we consider the fractional Poisson Process: this means that $K$ is
taken to be the  kernel
 $K_H$ associated to fractional Brownian motion (see \cite{du}):
\begin{displaymath}
  K_H(t,s)=\frac{1}{\Gamma(H+1/2)}(t-s)^{H-1/2}F(H-1/2,1/2-H,H+1/2,1-t/s),
\end{displaymath}
where $F(a,b,c,z)=\frac{\Gamma(c)}{\Gamma(b)\Gamma(c-b)}\int_0^1u^{b-1}(1-u)^{c-b-1}(1-zu)^{-a}\ 
du,$
 and
that $N$ is a Poisson process with constant intensity $\lambda.$ Then,  we know (see Decreusefond Üstünel \cite{du}) that
$\phi$ is given by: 
\begin{equation}\label{eq:3}
s\to\phi(s)=\frac{\Gamma(\frac{3}{2}-H)}{\Gamma(2-2H)} \frac{s^{\frac{1}{2}-H}}{\lambda}
\end{equation}
The theorem can be applied because $\phi \in \mathcal{L}^1(\nu)$ as
soon as $\frac{1}{2} < H < 1$ and in this case so does $\ln(1+\theta\phi).$
\end{remar}
We now prove that for $\theta >0,$  $\hat{\theta}_t$ is strongly consistant, i.e., that
$\hat{\theta}_t$ converges a.s. to $\theta.$
\begin{lemma}
  The process $\hat{\theta}_t$ is decreasing on any interval $[T_n,\, T_{n+1}[.$
\end{lemma}
\begin{proof}
  Let $\psi(x,y)=x/(1+xy),$ it is clear that for $x>0,$ the partial
  map $(y\mapsto \psi(x,y))$ is decreasing. We know that
  $\hat{\theta}_t$ is the solution of the equation
  \begin{equation}
    \sum_{T_j\le t}\psi(\phi(T_j),\hat{\theta}_t)=\int_0^t
    \phi(s)\lambda(s)\, ds.\label{eq:2}
  \end{equation}
For $t\in [T_n,T_{n+1}[,$ the number of terms in the left-hand-side in
\eqref{eq:2} is 
constant and the right-hand-side is an increasing function of $t.$ It
follows that $\hat{\theta}_t$ must decrease between $T_n$ and $T_{n+1}.$
\end{proof}
To go further, we need additional hypothesis on $\phi.$ Note that all
these hypothesis are satisfied by $\phi$ as defined in \eqref{eq:3}. 
We don't know whether $\hat{\theta}_{T_n}$ is decreasing sequence and
we cannot thus conclude  to the convergence of $\hat{\theta}_t.$
However, we now prove that it is a bounded process. 
\begin{lemma}
  Assume  $\theta >0$ and that $\phi$ and $\ln(1+\theta\phi)$ belong to
  $L^1([0,T],\lambda(s)ds)$ for any $T>0.$  Assume that 
  \begin{displaymath}
 \lim_{t\to \infty}\int_0^t \phi^2(s)\lambda(s)\, ds=\infty 
  \end{displaymath}
and that 
\begin{equation}\label{eq:13}
  \lim_{t\to \infty}\dfrac{\int_0^t \phi^{2+j}(s)\lambda(s)\, ds}{
    \int_0^t \phi^2(s)\lambda(s)\, ds}=0, \text{ for any } j>0.
\end{equation}
Then, $\{\hat{\theta}_t, \, t\ge 0\}$ is $\P_\theta$-a.s. bounded.
\end{lemma}
\begin{proof}
Let $M > \theta$ and consider 
$$
A_M=\{ \omega \in \Omega,\, \limsup_{n\to\infty}
\hat{\theta}_t(\omega) \geq M \}.
$$
On $A_M$, there exists  a sequence $\{ t_n, \, n \ge 1 \}$ of positive
reals such that $\hat{\theta}_{t_n} \geq M$ for any $n \ge 1$. Thus
\begin{equation*}
 \sum_{T_j \le t_n}\psi(\phi(T_j),\hat{\theta}_{t_n}) \leq 
 \sum_{T_j \le t_n}\psi(\phi(T_j),M)
\end{equation*}
By the very definition of $\hat{\theta}_t$, 
\begin{equation*}
    \sum_{T_j\le t_n}\psi(\phi(T_j),\hat{\theta}_{t_n})=\int_0^{t_n}
    \phi(s)\lambda(s)\, ds,
\end{equation*}
thus
\begin{multline}  \label{eq:11}
\int_0^{t_n} \phi(s)\lambda(s)\, ds- \int_0^{t_n}
\psi (\phi(s),M) (1+ \theta \phi(s)) \lambda(s) \, ds \\
\leq 
 \sum_{T_j \le t_n}\psi(\phi(T_j),M) -
\int_0^{t_n} \psi(\phi(s),M) (1+ \theta\phi(s))
\lambda(s)\, ds.
\end{multline}
Left-hand-side of Equation \eqref{eq:11} can be
simplified as
\begin{multline*}
\int_0^{t_n} \phi(s)\lambda(s)\, ds- \int_0^{t_n}
\psi(\phi(s),M) (1+ \theta \phi(s)) \lambda(s) \, ds \\
=
(M-\theta) \, \int_0^{t_n} \frac{\phi^2(s)\lambda(s)}{1+M\phi(s)} \,
ds.   
\end{multline*}
Thus, Equation \eqref{eq:11} now reads as 
\begin{multline}
  \label{eq:9}
  (M-\theta) \, \int_0^{t_n} \frac{\phi^2(s)\lambda(s)}{1+M\phi(s)} \,
  ds \le \sum_{T_j \le t_n}\psi(\phi(T_j),M)\\ -
\int_0^{t_n} \psi(\phi(s),M) (1+ \theta\phi(s))
\lambda(s)\, ds.
\end{multline}
Furthermore, since $\phi \ge 0$ and $M>\theta,$ we have
\begin{equation}
\label{eq:12}
\int_0^t \frac{\phi^2(s)\lambda(s)(1+\theta\phi(s))}{(1+M\phi(s))^2}
\, ds \le \int_0^t \frac{\phi^2(s)\lambda(s)}{1+M\phi(s)} \, ds.   
\end{equation}
On the other hand, Central Limit theorem for martingale, says that, $\P_\theta$ a.s.,
\begin{multline*}
  \dfrac{\sum_{T_j \le t_n} \psi(\phi(T_j),M) -
\int_0^{t_n}  \psi(\phi(s),M) (1+ \theta\phi(s))
\lambda(s)\, ds}{ \int_0^{t_n} \psi(\phi(s),M)^2  \lambda(s)(1+\theta \phi(s)) \, ds}
 \xrightarrow{n\to\infty}{}0,
\end{multline*}
Then, divide both sides of \eqref{eq:9} by $ \int_0^{t_n} \psi(\phi(s),M) \phi(s)\lambda(s) \,
ds$ and let $n$ go to infinity. In virtue of \eqref{eq:13} and \eqref{eq:12}, this yields to the
contradiction that on $A_M,$ $\theta \ge M.$ Finally, $\P_\theta(A_M)=0$
and thus $\{ \hat{\theta}_t,\, t > 0 \}$ is $\P_\theta$-a.s. bounded.
\end{proof}
\begin{theo}
Assume  $\theta >0$ and that $\phi$ and $\ln(1+\theta\phi)$ belong to
  $L^1([0,T],\lambda(s)ds)$ for any $T>0.$  Assume that 
  \begin{displaymath}
   \lim_{t\to \infty}\int_0^t \phi^2(s)\lambda(s)\, ds=\infty 
  \end{displaymath}
and that 
\begin{displaymath}
  \lim_{t\to \infty}\dfrac{\int_0^t \phi^{2+j}(s)\lambda(s)\, ds}{
    \int_0^t \phi^2(s)\lambda(s)\, ds}=0, \text{ for any } j>0.
\end{displaymath}
Then $\hat{\theta}_t$ tends $\text{P}_\theta$-a.s. to $\theta.$
\end{theo}
\begin{proof}
Let $N_t(\phi)=\sum_{T_n\le t}\phi(T_n).$
  According to \eqref{eq:2}, we have
  \begin{equation*}
    \begin{split}
0&  = \sum_{T_n\le t}\frac{\phi(T_n)}{1+\hat{\theta}_t\phi(T_n)}-\int_0^t \phi(s)\lambda(s)\, ds\\
& = N_t(\phi)-\int_0^t \phi(s)\lambda(s)(1+\theta\phi(s))\, ds\\
& - \hat{\theta}_t \Bigl(N_t(\phi^2)-\int_0^t
\phi^2(s)\lambda(s)(1+\theta\phi(s))\, ds\Bigr)\\
& +\hat{\theta}_t^2\sum_{T_n\le
  t}\frac{\phi(T_n)^3}{1+\hat{\theta}_t\phi(T_n)}-\theta\hat{\theta}_t\int_0^t
\phi^3(s)\lambda(s)\, ds\\
&+(\theta-\hat{\theta}_t)\int_0^t
\phi^2(s)\lambda(s)\, ds.      
    \end{split}
  \end{equation*}
After some simple algebra, this yields to 
\begin{multline}
  \label{eq:4}
  (\theta-\hat{\theta}_t)\bigl(\int_0^t
\phi^2(s)\lambda(s)\, ds +\hat{\theta}_t \int_0^t
\phi^3(s)\lambda(s)\, ds\bigr)\\
  \begin{aligned}
   & 
=-(N_t(\phi)-\int_0^t \phi(s)\lambda(s)(1+\theta\phi(s))\, ds)\\
&+  \hat{\theta}_t \Bigl(N_t(\phi^2)-\int_0^t
\phi^2(s)\lambda(s)(1+\theta\phi(s))\, ds\Bigr)\\
&-\hat{\theta}_t^2\Bigl(N_t(\phi^3)-\int_0^t
\phi^3(s)\lambda(s)(1+\theta\phi(s))\, ds\Bigr)\\
&-\hat{\theta}_t^2 \sum_{T_n\le
  t}\Bigl(\frac{\phi(T_n)^3}{1+\hat{\theta}_t\phi(T_n)}-\phi(T_n)^3\Bigr)+\theta \hat{\theta}_t^2 \int_0^t
\phi^4(s)\lambda(s)\, ds.
  \end{aligned}
\end{multline}
For any deterministic $\zeta,$ the process $\{N_t(\zeta)-\int_0^t
\zeta(s)\lambda(s)(1+\theta\phi(s))\, ds,\, t\ge 0\}$ is a martingale
 whose square bracket given by 
\begin{displaymath}
  \int_0^t \zeta(s)^2\lambda(s)(1+\theta\phi(s))\, ds.
\end{displaymath}
It follows from  hypothesis and Central Limit Theorem for
martingales that 
\begin{equation}
\label{eq:6}    \lim_{t\to \infty}\frac{N_t(\phi^j)-\int_0^t \phi(s)^j\lambda(s)(1+\theta\phi(s))\, ds}{ \int_0^t \phi^2(s)\lambda(s)\, ds}=0,
\end{equation}
for any $j\ge 1.$
Moreover, 
\begin{equation}
\label{eq:7}  \sum_{T_n\le
  t}\Bigl(\frac{\phi(T_n)^3}{1+\hat{\theta}_t\phi(T_n)}-\phi(T_n)^3\Bigr)
=\hat{\theta}_t \sum_{T_n\le
  t}\frac{\phi(T_n)^4}{1+\hat{\theta}_t\phi(T_n)}.
\end{equation}
Thus,
\begin{equation*}
  0\le \sum_{T_n\le
  t}\Bigl(\frac{\phi(T_n)^3}{1+\hat{\theta}_t\phi(T_n)}-\phi(T_n)^3\Bigr) \le \hat{\theta}_t N_t(\phi^4).
\end{equation*}
According to \eqref{eq:6} and to hypothesis \eqref{eq:13}, we have
\begin{equation}
  \label{eq:10}
  \lim_{t\to\infty}\frac{N_t(\phi^4)}{\int_0^t \phi^2(s)\lambda(s)\,
    ds}=0,\ \P_\theta \text{ a.s..}
\end{equation}
Divide \eqref{eq:4} by $\int_0^t \phi^2(s)\lambda(s)\, ds$ and
let $t$ goes to infinity, it follows from \eqref{eq:6} and
\eqref{eq:10} that $\hat{\theta}_t$ converges to $\theta.$   
\end{proof}
\def\polhk#1{\setbox0=\hbox{#1}{\ooalign{\hidewidth
  \lower1.5ex\hbox{`}\hidewidth\crcr\unhbox0}}}

\end{document}